\documentclass[12pt]{amsart}

\newtheorem{theorem}{Theorem}

\newtheorem{remark}[theorem]{Remark}
\newtheorem*{theorem*}{Theorem}
\newtheorem{proposition*}[theorem]{Proposition}

\newtheorem{example}[theorem]{Example}
\newtheorem{corollary}[theorem]{Corollary}

\theoremstyle{remark}

\usepackage{amssymb}

\def\CC{\mathbb C}

\def\NN{\mathbb N}
\def\ZZ{\mathbb Z}
\def\RR{\mathbb R}
\def\CC{\mathbb C}
\def\DD{\mathbb D}
\def\too{\longrightarrow}
\def\CO{\mathcal O}

\def\bs{\boldsymbol}
\def\Sp{|\bs{p}|}
\def\Sq{|\bs{q}|}
\def\phi{\varphi}
\def\dist{\operatorname{dist}}
\def\Ree{\operatorname{Re}}

\def\inte{\operatorname{int}}

\title[]
{The multipole Lempert function is monotone under inclusion
of pole sets}
\author{Nikolai Nikolov and Peter Pflug }

\address
{Institute of Mathematics and Informatics\\ Bulgarian
Academy of Sciences\\
1113 Sofia, Bulgaria} \email{nik@math.bas.bg}
\address{Carl von Ossietzky Universit\"at Oldenburg\\
Fachbereich Mathematik\\ Postfach 2503\\ D-26111 Oldenburg,
Germany}
\email{pflug@mathematik.uni-oldenburg.de}

\begin{document}


\begin{abstract}
We prove that the multipole Lempert function is monotone
under inclusion of
pole sets.
\end{abstract}

\maketitle

Let $D$ be a domain in $\CC^n$ and let $A=(a_j)_{j=1}^l$,
$1\le
l\le\infty$, be a countable (i.e. $l=\infty$) or a
non--empty
finite subset of $D$ (i.e. $l\in\NN$). Moreover, fix a
function
$\bs p:D\too\RR_+$ with
$$
\Sp:=\{a\in D:\bs p(a)>0\}=A.
$$
$\bs p$ is called a {\it pole function for $A$} on $D$ and
$\Sp$
its {\it pole set}. In case that $B\subset A$ is a
non--empty
subset we put $\bs p_B:=\bs p$ on $B$ and $\bs p_B:=0$ on
$D\setminus B$. $\bs p_B$ is a pole function for $B$.

For $z\in D$ we set
$$
l_D(\bs p,z)=\inf\{\prod_{j=1}^l|\lambda_j|^{\bs p(a_j)}\},
$$
where the infimum is taken over all subsets
$(\lambda_j)_{j=1}^l$
of $\DD$ (in this paper, $\DD$ is the open unit disc in
$\CC$) for
which there is an analytic disc $\phi\in\CO(\DD,D)$ with
$\phi(0)=z$ and $\phi(\lambda_j)=a_j$ for all $j$. Here we
call
$l_D(\bs p,\cdot)$ {\it the Lempert function with $\bs
p$-weighted
poles at $A$} \cite{Wik1,Wik2}; see also \cite{jarpfl},
where this
function is called  {\it the Coman function for $\bs p$}.

Recently, F.~Wikstr\"om \cite{Wik1} proved that if $A$ and
$B$ are finite
subsets of a convex domain $D\subset\CC^n$ with
$\varnothing\neq B\subset A$
and if $\bs p$ is a pole function for $A$, then $l_D(\bs
p,\cdot)\le l_D(\bs
p_B,\cdot)$ on $D$.

On the other hand, in \cite{Wik2} F.~Wikstr\"om gave an
example of
a complex space for which this inequality fails to hold and
he
asked whether it remains to be true for an arbitrary domain
in
$\CC^n$.

The main purpose of this note is to present a positive
answer to that question, even
for countable pole sets (in particular, it follows that the
infimum in the
definition of the Lempert function is always taken over a
non-empty set).

\begin{theorem}  For any domain $D\subset\CC^n$, any
countable or non-empty finite subset $A$
of $D$, and any pole function $\bs p$ for $A$ we have
$$
l_D(\bs p,\cdot)=\inf\{l_D(\bs p_B,\cdot):\varnothing\neq B
\text{ a
finite subset of } A\}.
$$

Therefore, $$ l_D(\bs p,\cdot)=\inf\{l_D(\bs
p_B,\cdot):\varnothing\neq
B\subset A\}. $$
\end{theorem}

\vskip 0.5cm

The proof of this result will be based on the following

\

\begin{theorem*}[Arakelian's Theorem \cite{Ara}]  Let
$E\subset\Omega\subset\CC$ be a relatively
closed subset of the domain $\Omega$. Assume that
$\Omega^*\setminus E$ is connected and
locally connected. (Here $\Omega^*$ denotes the one--point
compactification of $\Omega$.)

 If $f$ is a complex-valued continuous function on $E$ that
is
holomorphic in the interior of $E$ and if $\varepsilon>0$,
then there is
a $g\in\CO(\Omega)$ with $|g(z)-f(z)|<\varepsilon$ for any
$z\in E$.
\end{theorem*}
\

\begin{proof} Fix a point $z\in D$. First, we shall verify
the inequality
$$
l_D(\bs p,z)\le\inf\{l_D(\bs p,z):\varnothing\neq B \text{
a
finite subset of } A\}:\leqno{(1)}
$$

Take a non--empty proper finite subset $B$ of $A$. Without
loss of
generality we may assume that $B=A_m:=(a_j)_{j=1}^m$ for a
certain
$m\in\NN$, where $A=(a_j)_{j=1}^l$, $m<l\le\infty$.

Now, let $\phi:\Bbb D\to D$ be an analytic disc with
$\phi(\lambda_j)=a_j,$ $1\le j\le m,$ where $\lambda_0:=0$ and
$a_0:=z.$ Fix $t\in[\max_{0\le j\le m}|\lambda_j|,1)$ and put
$\lambda_j=1-\frac{1-t}{j^2}$,\; $j\in A(m)$, where
$A(m):=\{m+1,\dots,l\}$ if $l<\infty$, respectively
$A(m):=\{j\in\NN:j>m\}$ if $l=\infty$. Consider a continuous curve
$\phi_1:[t,1)\to D$ such that $\phi_1(t)=\phi(t)$ and
$\phi_1(\lambda_j)=a_j$, $j\in A(m)$. Define
$$
f=\begin{cases}
\phi|_{\overline{t\Bbb D}}\\
\phi_1|_{[t,1)}
\end{cases}
$$
on the set $F_t=\overline{t\Bbb D}\cup[t,1)\subset\DD$.
Observe that
$F_t$ satisfies the geometric condition in Arakelian's
Theorem.

Since $(\lambda_j)_{j=0}^l$ satisfy the Blaschke condition,
for any $k$ we
may find a Blaschke product $B_k$ with zero set
$(\lambda_j)_{j=0,j\neq
k}^l$. Moreover, we denote by $d$ the function
$\dist(\partial D,f)$ on
$F_t$, where the distance arises from the $l^\infty$--norm.
Let $\eta_1$,
$\eta_2$ be  continuous real-valued functions on $F_t$ with
\begin{gather*}
\eta_1,\eta_2\le\log\frac{d}{9},\qquad
\eta_1,\eta_2=\min_{\overline{t\DD}}
\log\frac{d}{9}\hbox{ on }(\overline{t\DD}),\\
\text{ and }\quad
\eta_1(\lambda_j)=\eta_2(\lambda_j)+\log(2^{-j-1}|B_j(\lambda_j)|),\quad
j\in A(m).
\end{gather*}

Applying three times
Arakelian's theorem, we may find functions $\zeta_1,
\zeta_2\in\CO(\DD)$ and a holomorphic mapping $h$ on $\DD$
such that
$$
|\zeta_1-\eta_1|\le 1,\;|\zeta_2-\eta_2|\le 1,\hbox{ and
}|h-f|\le
\varepsilon|e^{\zeta_1-1}|\le \varepsilon e^{\eta_1}\text{
on } F_t,
$$
where $\varepsilon:=
\min\{\tfrac{|B_j(\lambda_j)|}{2^{j+1}}:j=0,\dots,m\}<1$
(in the last case
apply Arakelian's theorem componentwise to the
mapping $e^{1-\zeta_1}f$).

In particular, we have
\begin{gather*}
|h-f|\le\frac{d}{9}\quad \text{ on } F_t,\\
|\gamma_j|\leq
e^{\eta_1(\lambda_j)}2^{-j-1}|B_j(\lambda_j)|
=e^{\eta_2(\lambda_j)}2^{-j-1}|B_j(\lambda_j)|,\quad
j=0,\dots,m,\\
|\gamma_j|\le e^{\eta_1(\lambda_j)}=
e^{\eta_2(\lambda_j)}2^{-j-1}|B_j(\lambda_j)|
,\;j\in A(m),
\end{gather*}
where $\gamma_j:=h(\lambda_j)-f(\lambda_j)$, $j\in\ZZ_+$ if
$l=\infty$, respectively $0\leq j\leq l$ if $l\in\NN$.
Then, in virtue of
$e^{\eta_2(\lambda_j)}
\le e^{1+\Ree \zeta_2(\lambda_j)}$, the function
$$
g:=e^{\zeta_2}\sum_{j=0}^l\frac{B_j
}{e^{\zeta_2(\lambda_j)}B_j(\lambda_j)}\gamma_j
$$
is
holomorphic on $\DD$ with $g(\lambda_j)=\gamma_j$ and $$
|g|\le e^{\Ree\zeta_2+1}\le
e^{\eta_2+2}\le\frac{e^2}{9}d\quad
\text{ on } F_t. $$ For $q_t:=h-g$ it follows that
$q_t(\lambda_j)=f(\lambda_j)$ and $$
|q_t-f|\le\frac{e^2+1}{9}d<d\quad\text{ on }F_t. $$

Thus we have found a holomorphic
mapping $q_t$ on $\DD$ with $q_t(\lambda_j)=a_j$ and
$q_t(F_t)\subset D$.
Hence there is a simply connected domain
$E_t$ such that $F_t\subset E_t\subset\DD$ and
$q_t(E_t)\subset D$.

Let $\rho_t:\DD\to E_t$ be the Riemann mapping with
$\rho_t(0)=0,$ $\rho_t'(0)>0$ and
$\rho_t(\lambda_j^t)=\lambda_j$. Considering the analytic
disc
$q_t\circ\rho_t:\DD\to D$, we get that
$$
l_D(\bs p,z)\le\prod_{j=1}^l|\lambda_j^t|^{\bs
p(a_j)}\leq\prod_{j=1}^m|\lambda_j^t|^{\bs p(a_j)}.
$$
Note that by the Carath\'eodory kernel theorem, $\rho_t$
tends, locally uniformly, to the
identity map of $\DD$ as $t\to 1$. This shows that the last
product converges to
$\prod_{j=1}^m|\lambda_j|^{\bs p(a_j)}$. Since $\phi$ was
an arbitrary
competitor for $l_G(\bs p_{A_m},z)$, the inequality (1)
follows.

On the other hand, the existence of an analytic disc whose
graph contains
$A$ and $z$ implies that
$$
l_D(\bs p,z)\ge\limsup_{m\to\infty}l_D(\bs p_{A_m},z),
$$
which completes the proof.
\end{proof}

\begin{remark}{\rm Looking at the above proof shows that we
have proved an
approximation and simultaneous interpolation result, i.e.
the
constructed function $q_t$ approximates and interpolates
the given
function $f$.

We first mention that the proof of Theorem 1 could be
simplified
using a non-trivial result on interpolation sequences due
to
L.~Carleson (see, for example, Chapter 7, Theorem 3.1 in
\cite{and}).

Moreover, it is possible to prove the following general
result
which extends an approximation and simultaneous
interpolation
result by P.~M.~Gauthier, W.~Hengartner, and
A.~A.~Nersesyan (see \cite{gauhen}, \cite{ner}):

Let $D\subset\CC$ be a domain, $E\subset D$ a relatively
closed
subset satisfying the condition in Arakelian's theorem,
$\Lambda\subset E$  such that $\Lambda$
has no accumulation point in $D$ and
$\Lambda\cap\inte E$ is a finite set. Then for given
functions $f,h\in\mathcal C(E)\cap\CO(\inte E)$ there
exists a
$g\in\CO(D)$ such that
$$
|g-f|<e^{\Ree h} \text { on } E,\quad
 \text{ and } \quad
f(\lambda)=g(\lambda),\;\lambda\in\Lambda.
$$
It is even possible to prescribe a finite number of the
derivatives for $g$ at all the points in $\Lambda$. }

\end{remark}

\
As a byproduct we get the following result.

\begin{corollary} Let $D\subset\CC^n$ be a domain and let
$\bs p, \bs q:D\to\RR_+$ be two
pole functions on $D$ with $\bs p\leq\bs q$, $\Sq$ at most
countable. Then
$l_D(\bs q,z)\leq l_D(\bs p,z)$, $z\in D$.
\end{corollary}

Hence, the Lempert function is monotone with respect to
pole functions with an
at most countable support.

\begin{remark}{\rm The Lempert function is, in general, non
strictly monotone under
inclusion of pole sets. In fact, take $D:=\DD\times\DD$,
$A:=\{a_1,a_2\}\subset\DD$, where $a_1\neq a_2$, and
observe that
$l_D(\bs p,(0,0))=|a_1|$, where $\bs
p:=\chi|_{A\times\{a_1\}}$
(use the product property in \cite{dietra}, \cite{jarpfl},
\cite{nikzwo}).
}\end{remark}

\begin{remark}{\rm Let $D\subset\CC^n$ be a domain, $z\in
D$, and let now $\bs
p:D\too\RR_+$ be a "general" pole function, i.e. $\Sp$ is
uncountable. Then there are two cases:

1) There is a $\phi\in\CO(\DD,D)$ with
$\phi(\lambda_{\phi,a})=a$, $\lambda_{\phi,a}\in\DD$,
for all $a\in\Sp$ and $\phi(0)=z$. Defining
\begin{multline*}
l_D(\bs
p,z):=\inf\{\prod|\lambda_{\psi,a}|^{\bs
p(a)}:\psi\in\CO(\DD,D) \text{ with
}\\ \psi(\lambda_{\psi,a})=a \text{ for all } a\in\Sp,
\psi(0)=z\}
\end{multline*}
then $l_D(\bs p,z)=0$.

Observe that $l_D(\bs p,z)=\inf\{l_D(\bs p_B,z):\varnothing
\neq B\text{ a finite
subset of } A\}$.

2) There is no analytic disc as in 1). In that case we may
define
$^{1)}$\footnote{1) Compare the definition of the Coman
function (for the second
case) in \cite{jarpfl}.}
$$
l_D(\bs p,z):=\inf\{l_D(\bs p_B,z):\varnothing\neq B\text {
a finite subset of
} A\}.
$$
Example \ref{ex} below may show that the definition in 2)
is more sensitive than
the one used in \cite{jarpfl}.
}
\end{remark}

\vskip 0.5cm

Before giving the example we use the above definition of
$l_D(\bs p,\cdot)$
for an arbitrary pole function $\bs p$ to conclude.

\begin{corollary} Let $D\subset\CC^n$ be a domain and let
$\bs p, \bs q :D\too\RR_+$
be arbitrary pole functions with $\bs p\leq\bs q$. Then
$l_D(\bs q,\cdot)\leq
l_D(\bs p,\cdot)$.
\end{corollary}

\begin{example}\label{ex} {\rm Put $D:=\DD\times\DD$ and
let $A\subset\DD$ be uncountable, e.g. $A=\DD$.
Then there is no $\phi\in\CO(\DD,D)$ passing through
$A\times\{0\}$ and
$(0,w)$, $w\in\DD_*$. Put $\bs p:=\chi|_{A\times\{0\}}$ on
$D$ as a pole function.

Let $B\subset A$ be a non--empty finite subset. Then
applying the
product property (see \cite{dietra}, \cite{jarpfl},
\cite{nikzwo}), we get
$$
l_D(\bs
p_B,(0,w))=g_D(B\times\{0\},(0,w))=\max\{g_{\DD}(B,0),g_{\DD}(0,w)\},
$$
where $g_{\DD}(A,\cdot)$ denotes the Green function in
$\DD$ with respect to
the pole set $A$.

Therefore, $l_D(\bs p,(0,w))=|w|$ $^{1)}$.}
\end{example}

To conclude this note we mention that the Lempert function
is not
holomorphically contractible, even if the holomorphic map
is a proper
covering.

\begin{example} {\rm Let $\pi:\DD_*\too\DD_*$,
$\pi(z):=z^2$. Obviously, $\pi$ is
proper and a covering. Fix two different points $a_1,
a_2\in\DD_*$ with $a_1^2=a_2^2=:c$. For a point
$z\in\DD_*$ we know that
$$
l_{\DD_*}(\chi|_{\{c\}},z^2)=\min\{l_{\DD_*}(\chi|_{\{a_1\}},z),l_{\DD_*}(\chi|_{\{a_2\}},z)\}
\geq l_{\DD_*}(\chi|_A,z),
$$
where $A:=\{a_j:j=1,2\}$ and $\chi|_B$  is the
characteristic function for the set $B\subset D_*$.
Assume that $l_{\DD_*}(\chi|_{\{a_1\}},z)\leq
l_{\DD_*}(\chi|_{\{a_2\}},z)$.
Then this left side is nothing than the classical Lempert
function for the pair $(a_1,z)$.
Recalling how to calculate it via a covering map one easily
concludes that $l_{\DD_*}(\chi|_A,z)<
l_{\DD_*}(\chi|_{\{a_1\}},z)$. Hence $l_{\DD_*}(\bs p
,\pi(z))>l_{\DD_*}(\bs p\circ\pi,z)$, where $\bs
p:=\chi|_{\{c\}}$.

Therefore, the Lempert function with a multipole behaves
worse than the Green
function with a multipole.}
\end{example}

{\bf Acknowledgments.} This note was written during the
stay of
the first named author at the University of Oldenburg
supported by a
grant from the DFG (September - October 2004). He likes to
thank both
institutions for their support.

We thank the referee for pointing out an error in the
first version of this paper.

\end{document}